\documentclass[a4paper,12pt]{amsart}
\usepackage{hyperref}
\usepackage{bm}
\usepackage{mathtools}
\usepackage{graphicx}
\usepackage{amsmath,amssymb,amsthm,amsfonts}
\usepackage{amssymb}
\usepackage[active]{srcltx}
\usepackage{bbold}
\usepackage{color}
\newtheorem{thm}{Theorem}[section]

\newtheorem{rem}[thm]{Remark}

\newtheorem{Oquestion}{Question}
\newtheorem{OAnswer}{Answer}

\usepackage{tikz}
\usepackage{standalone}
\usepackage{graphicx}
\usepackage{subcaption}
\usepackage{multicol,float}

\numberwithin{equation}{section}

\newcommand{\R}{\mathbb R}

\begin{document}
	
	\parindent 0pc
	\parskip 6pt
	\overfullrule=0pt
	
	\title{\quad The Eleven Questions : March 2016 - October 2026}
	
	\author{Alberto Farina$^{\ast}$}
	\address{$^{\ast}$  Universit\'e de Picardie Jules Verne, LAMFA, CNRS UMR 7352, 33, rue Saint-Leu 80039 Amiens, France}\email{alberto.farina@u-picardie.fr}

	
	\date{}
	

	\begin{abstract}
We review the eleven open questions in Appendix H of Farina’s manuscript on entire solutions
of the Allen–Cahn equation \cite{Far-9/3/16}.
	\end{abstract}

	\maketitle

\bigskip

Farina’s manuscript \cite{Far-9/3/16} ends with eleven open questions concerning entire solutions to the Allen-Cahn equation
\begin{equation}\label{equaz-Allen-Cahn}
-\Delta u = u - u^3  \qquad \text{in} \qquad \R^N , \qquad N \ge 1.
\end{equation}
In this Note we use the notation, terminology and definitions introduced in \cite{Far-9/3/16}. 
We follow the order in which the open problems are presented in Appendix H of \cite{Far-9/3/16}, as well as their original formulation.

We will provide a summary of all eleven questions, distinguishing between
published results between March 2016 and October 2026, new claims contained in recent preprints, and inferences that combine them with the results demonstrated in the manuscript \cite{Far-9/3/16}.

Finally, I would like to point out that while drafting this Note, we chose not to take into account the sources \cite{OpenAI-1, OpenAI-2}, nor the very recent preprints that are based on the results and tools contained in those as-yet-unverified sources.

\section{Open Question 1}

\begin{Oquestion} Assume $ N\ge 9$. In view of the gap phenomena for local minimizers (cf. Theorem 2.3 of \cite{Far-9/3/16}) and the local minimizer built by M. del Pino, M. Kowalczyk and J. Wei in \cite{DKW}, it is easily seen that the \textsf{second energy level for local minimizers} is well defined. As for $ \ell_2$, the second positive energy level for finite energy solutions, it is defined by \\
\begin{equation}\label{defLiv2-min-OQ}
{\ell}_2(\mathcal{M}) : = \inf \left\lbrace \, L_u > \ell_{min} \,\, : \,\, u \in \mathcal{M}  \, \right\rbrace \footnote{Recall that $\mathcal{M}$ denotes the set of nonconstant local minimizers.}
\end{equation} 

\medskip

Is ${\ell}_{2}(\mathcal{M})$ achieved by a member of $\mathcal{M}$ when $N \ge 9 $ ? In the affirmative, compute the exact value.
\end{Oquestion}
\noindent A possible candidate for $N=9$, is the local minimizer built in \cite{DKW}. 

\smallskip

\begin{OAnswer} Attainment and the exact second level energy in the class of local minimizers remain open in general for $N \geq9$. 

\end{OAnswer}

\section{Open Question 2}

\begin{Oquestion} Does there exist ${\ell}_{2}(\mathcal{M})$ in $ \R^8$ ? 
\end{Oquestion}

\begin{OAnswer} Yes, the second level exists in $\R^8$ and it is attained. 
\end{OAnswer}
\noindent \textit{Existence - } In 2017, Liu, Wang and Wei \cite{LWW-2017} constructed non-planar local minimizers for any $ N \geq 8$. More precisely, given any strict area-minimizing Lawson’s cone, Liu, Wang and Wei proved the existence of a family of local minimizers whose nodal sets are asymptotic to this cone. 
From this, we infer the existence of ${\ell}_{2}(\mathcal{M})$.

\noindent \textit{Attainment of the second energy level - } We deduce this result by combining Theorem 5.1 of \cite{Edelen}, Theorem 2.1 of \cite{Far-9/3/16}, the classical regularity theory for area-minimizing cones and the existence of the second energy level for local minimizers just established. 
Edelen’s Theorem 5.1 of \cite{Edelen} gives local finiteness of densities of stable minimal cones in $\R^8$ that are smooth away from their vertex and with unit multiplicity. Every area-minimizing cone in $\R^8$ belongs to this class. Items $i),ii)$ and $iv)$ of Theorem 2.1 of \cite{Far-9/3/16} imply that any blow-down of a local minimizer $u$ is an area-minimizing cone $C$, with vertex at the origin and 
satisfying  the universal bound
$$
1 \leq \Theta(C,0) = \frac{L_u(8)}{\ell_{min}(8)} \leq \frac{8\omega_8}{2\omega_7} = \frac{35 \pi}{32} .
$$  
It follows that the realised positive energy spectrum of $\mathcal{M}$ in dimension eight is finite. This fact, together with the existence of the second energy level for local minimizers just established, proves that ${\ell}_{2}(\mathcal{M})$ is attained in dimension $8$.

\noindent The exact value of ${\ell}_{2}(\mathcal{M})$ has yet to be determined. 

\begin{rem}\label{Rem-OQ2}
The previous argument also ensures that the number of energy levels in the class $ \mathcal{M}$ of local minimizers is finite for $N=8$.  
\end{rem}

\section{Open Question 3}

\begin{Oquestion} Does there exist a local minimizer in $\R^8$ which is not 1D ?
\end{Oquestion} 

\noindent A possible candidate is the saddle-type solution built by X. Cabr\'e and J. Terra in \cite{CT1} and studied in detail in Appendix D  of \cite{Far-9/3/16}. 

\begin{OAnswer} Yes, there exists a local minimizer in $\R^8$ which is not planar.
\end{OAnswer} 
In 2017, Liu, Wang and Wei \cite{LWW-2017} constructed a non-planar 
$O(4) \times O(4)$-invariant local minimizer asymptotic to the Simons cone 
$\mathcal{C}_{4,4}$. 

\noindent We point out that the $O(4) \times O(4)$-invariant local minimizer constructed in \cite{LWW-2017} may not be the saddle-type solution built in \cite{CT1} (see Remark 2 of  \cite{LWW-2017}). 

\section{Open Question 4}

\begin{Oquestion} Prove (or disprove) that there are only a finite number of energy levels in the set $\mathcal{M}$ of the local minimizers.
\end{Oquestion} 

\noindent For $ N\le 7$, there is only one energy level, the minimum one $\ell_{min}$, thanks to item iii) of Theorem 2.1 of \cite{Far-9/3/16}. The upper bound proved in item iv) of Theorem 2.1 of \cite{Far-9/3/16} suggests that this might be true also for $ N \ge 8$. 

\begin{OAnswer} The answer is affirmative for $N=8$. See Remark \ref{Rem-OQ2}.
\end{OAnswer}

\noindent The general question is still open for $N \geq 9$. 

\section{Open Question 5}

\begin{Oquestion} Prove (or disprove) that there are only a finite number of energy levels in the set $\mathcal{M}_m$ of the monotone local minimizers.
\end{Oquestion} 

\noindent For $ N\le 8$, thanks to item iii) of Theorem 2.1 of \cite{Far-9/3/16} and to Theorem 2.5 of \cite{Far-9/3/16}, there is only one energy level, the minimum one 
$\ell_{min}$.

\begin{OAnswer} For $N=9$ there are only a finite number of energy levels in the set $\mathcal{M}_m$. 
\end{OAnswer} 

We deduce this result from the results of the manuscript \cite{Far-9/3/16} and Theorem 5.1 of \cite{Edelen}. Indeed, Theorem 2.7 of the manuscript \cite{Far-9/3/16} reduces the blow-down of a non-planar monotone local minimizer $u$ in $\R^9$ to $C'\times \R$, where $C'$ is an area-minimizing cone in $\R^8$. Therefore, by Theorem 2.1 (items $(ii)$ and $(iv)$) of \cite{Far-9/3/16} we get
$$
\frac{9}{2} \frac{\omega_9}{\omega_8}\ell_{min}(9) \geq L_u = \Theta(C' \times \R,0) \ell_{min}(9) = \Theta(C',0') \ell_{min}(9)
$$   
and so, $\Theta(C',0') \leq \frac{9}{2} \frac{\omega_9}{\omega_8} = \frac{128}{35}$. 
But then, we also know \cite[Theorem 5.1]{{Edelen}} that, 
for stable minimal cones in $\R^8$ (with unit multiplicity) that are regular outside their vertex and whose densities lie below a fixed threshold, the set of values taken by the densities is finite. The desired conclusion then follows.  

\noindent The general finiteness problem remains open for $N \geq 10$.

\begin{rem}\label{Rem-OQ5}
Since there exist non-planar monotone local minimizers in $\R^9$ (see the formulation of Open Problem 1), the second positive level is also reached within the class $\mathcal{M}_m$. This result must be distinguished from the unresolved attainment problem in the full class 
$\mathcal{M}$.
\end{rem}

\section{Open Question 6}

\begin{Oquestion} Do there exist nonconstant, non-monotone local minimizers for $N \ge 9$ ? 
\end{Oquestion} 

\noindent A possible candidate is the saddle-type solution built in \cite{CT1}, at least in even dimensions $N$.   

\begin{OAnswer} Yes for any $ N \geq 9$. 
\end{OAnswer} 

The non-planar local minimizers $u$ constructed in \cite{LWW-2017} for any $N \geq 9$ satisfy $u(x) = u(-x)$ for any $ x \in \R^N$, therefore they have a critical point at the origin. Hence, none of them is monotone.

\section{Open Question 7}

\begin{Oquestion} Assume $ N \ge 4$. Prove (or disprove) that every monotone solution has finite energy.
\end{Oquestion} 

\noindent As pointed out in Remark 15 of \cite{Far-9/3/16}, this is true for $ N\le 3$.

\begin{OAnswer} 
The answer is affirmative for $N=4$ on the basis of the statements in the two independent preprints \cite{CF-RFF-SS,LLWWWW} from September 2026.
\end{OAnswer} 
The two independent preprints \cite{CF-RFF-SS,LLWWWW} from September 2026 show that every stable solution in $\R^3$ is one-dimensional and use this result to establish the validity of De Giorgi’s conjecture for monotone solutions in $\R^4$. Consequently, the answer is affirmative in dimension $4$.

\section{Open Question 8}

\begin{Oquestion} Assume $ N \ge 2$ and let $u$ be a nonconstant finite-energy solution of the Allen-Cahn equation. Prove (or disprove) that
\[
u \, \textit{ is 1D} \quad \Longleftrightarrow \quad L_u = \ell_{min}. 
\]
\end{Oquestion} 

\noindent Note that this is true when $u$ is also a local minimizer, as proved in Theorem 2.2 of \cite{Far-9/3/16}. It is also true for stable solutions when $ N =2,3$ (cf. Theorem  6.1 of \cite{Far-9/3/16}). 

\begin{OAnswer} The answer is affirmative in any dimension. 
\end{OAnswer} 

K. Wang \cite{K-Wang-sav} proves that, for every $N \geq 2$, there is $\tau_N > 0$ such that any nonconstant entire solution $u$ satisfying
$$
L_u(N) \leq (1 + \tau_N) \ell_{min}(N)
$$
is of the form $u(x) = \tanh\big (\frac{(\nu\cdot x) +\alpha}{\sqrt 2} \big)$, where 
$\alpha \in \R$ and $ \nu $ is a unit vector of $\R^N$. The converse implication is well known.

\section{Open Question 9}

\begin{Oquestion} Assume $ N \ge 2$ and let $u$ be a nonconstant monotone solution of the Allen-Cahn equation. Prove (or disprove) that
\[
u \, \textit{ is 1D} \quad \Longleftrightarrow \quad L_u = \ell_{min}. 
\]
\end{Oquestion} 
\noindent Note that this is true when $u$ is also a local minimizer, as proved in Theorem 2.2 of \cite{Far-9/3/16}.

\begin{OAnswer} The answer is affirmative in any dimension. 
\end{OAnswer} 

See the answer to the Open Question 8. 

\section{Open Question 10}

\begin{Oquestion}Assume $N=2$. Do there exist entire solutions of the Allen-Cahn equation having a Pythagorean conical limit interface $V$ ? 
\end{Oquestion}
\begin{OAnswer} General existence of Pythagorean limits remains unresolved.
\end{OAnswer}

\section{Open Question 11}

\begin{Oquestion} 
\noindent Assume $ N\ge 5$. Is the second positive energy level ${\ell}_{2}$ achieved by a member of $\mathcal{F}_p$ ? \footnote{Recall that $\mathcal{F}_p$ is the class of solutions to \eqref{equaz-Allen-Cahn} satisfying $ 0 < L_u < +\infty $.} In the affirmative, compute the exact value. 
\end{Oquestion}

\noindent A good candidate is the saddle-type solution built in \cite{CT1}, at least in even dimensions $N.$  

\begin{OAnswer} Attainment and the exact value of the second energy level remain open in general for 
$N \geq5 $.
\end{OAnswer}

\vskip 1truecm

\section{Recent symmetry and rigidity results and their relation to 2016 manuscript}

\bigskip

\begin{enumerate}

\item[1)] Florit–Simon and Serra, in their preprint \cite[Theorem 1.3]{{FSF}}, classify stable solutions of finite energy in $\R^4$. These solutions are one-dimensional, that is, they are either planar heteroclines or the constants $ \pm 1$.
This strengthens Theorem 6.2 of \cite{Far-9/3/16}, which provides a planar blow-down support in dimension $N \leq 7$, but does not determine its multiplicity or classify the solution. 
Liu, Luo, Wang, Wei, Wei and Wu in the very recent preprint \cite[Theorem 1.2]{{LLWWWW-new}} extend the classification of stable solutions based on this energy hypothesis to all $N \leq 7$. 

\item[2)] Florit–Simon, in his preprint \cite[Theorem 1.1]{{FS}}, obtains the same classification result for finite energy solutions with finite Morse index in $\R^4$. The conditional extension to dimensions $ 4 \leq N \leq 7$ in \cite[Theorem 1.2]{{FS}}, combined with the classification of stable solutions of finite energy established in \cite{{LLWWWW-new}}, provides the classification for finite energy solutions with finite Morse index over the entire interval; the conclusion is also stated directly in \cite[Corollary 6.8]{{LLWWWW-new}}.
On the other hand, for $N=3$ solutions with finite energy and finite Morse index may be non-planar. The asymptotic structure of their nodal sets is described in \cite[Theorem 1.4]{{FS}}, which refines the description of the limit support contained in Theorem 7.1 of \cite{Far-9/3/16}.

\item[3)] For nonconstant finite-energy solutions with finite Morse index, Theorem 1.3 of \cite{{LLWWWW-new}} guarantees multiplicity equal to one for every varifold blow-down in the interval $4 \leq N \leq 10$. These results refine the results in Sections 6 and 7 of \cite{Far-9/3/16}. Observe that, multiplicity one neither gives a cone independent of the scale sequence nor proves minimality of the solution.

\item[4)] We deduce the \textit{validity of De Giorgi’s conjecture up to dimension $N=8$ for monotone solutions with finite energy}. This result extends Theorem 10.3 of  \cite{Far-9/3/16}, which held for $N=4$. 

Indeed, first recall that the limit profiles 
\begin{equation}\label{N-profiles-conj}
\begin{split}
&\underline u(x'):=\lim_{x_N\rightarrow-\infty} u(x',x_N), 
\\ 
&\overline u(x'):=\lim_{x_N\rightarrow+\infty} u(x',x_N),
\end{split}
\end{equation}
are stable solutions to the Allen-Cahn equation in $\R^{N-1}$ and then that, according to Theorems 1.5 and 1.6 of \cite{Far-9/3/16}, the two sequences of solutions defined by 
$u^{\pm}_{n}(x',x_N) := u(x',x_N \pm n)$, $n \geq1$, converge (up to a subsequence) to the solutions 
$u^+(x',x_N) := \overline u(x')$ and $u^-(x',x_N) := \underline u(x')$ in $C^2_{loc}(\R^N)$, with $L_{u^{\pm}} \leq L_u < \infty$. \\
Now we observe that, for $R>0$, 
\begin{align*}
&2R \int_{B^{N-1}(0',R)} \left( \frac{\vert \nabla \overline{\underline u} \vert^2}{2} + \frac{(1-\overline{\underline u}^2)^2}{4}\right)  dx' \\
& = \int_{B^{N-1}(0',R) \times (-R,R)} \left( \frac{\vert \nabla u^{\pm} \vert^2}{2} + \frac{(1-{u^{\pm}}^2)^2}{4}\right)  dx \\
& \leq \int_{B^{N}(0,\sqrt{2}R)} \left( \frac{\vert \nabla u^{\pm} \vert^2}{2} + 
\frac{(1-{u^{\pm}}^2)^2}{4}\right)  dx \leq L_{u^{\pm}} (\sqrt{2}R)^{N-1}  \leq L_u (\sqrt{2}R)^{N-1} 
\end{align*}
which in turn, immediately gives that the limit profiles $\underline u$ and 
$\overline u$ are finite-energy solutions in $\R^{N-1}$. Since $N-1 \leq 7$, they must therefore depend on at most one variable, according to the classification results presented
in item $1)$ above. The desired conclusion then follows by a direct application of the limit profile criterion of Farina and Valdinoci \cite[Theorem 1.1]{{Far-Val-TAMS}}.

\item[5)] We deduce that \textit{the set of possible energy levels attained in the class of stable solutions (resp. with finite Morse index) with finite energy in $\R^8$ is discrete \footnote{i.e., a subset of $\R$ consisting solely of isolated points.}; therefore, the second positive level is attained within each of the two classes under consideration.} 

Indeed, by Theorem 1.3 of \cite{LLWWWW-new} we know that every varifold blowdown of a nonconstant solution $u$ is a stable minimal cone $C$ with multiplicity one and, by item $iv)$ of Theorem 6.2 (resp. Theorem 7.1) of \cite{Far-9/3/16} we know that those cones are smooth outside their vertex. But then, we also know \cite[Theorem 5.1]{{Edelen}} that, 
for stable minimal cones in $\R^8$ (with unit multiplicity) that are regular outside their vertex, the set of values taken by the densities is discrete. The first desired conclusion then follows since $ L_u =  \Theta(C,0) {\ell_{min}(8)}$ by item $ii)$ of Theorem 6.2 (resp. Theorem 7.1) of \cite{Far-9/3/16}. Since nonplanar local minimisers in $\R^8$  belong to both classes and have energy above $\ell_{min}(8)$, the discreteness result just proven immediately implies that the second positive level is attained within each of the two classes under consideration.
\begin{rem} 
The argument presented above does not identify either of the two energy values, does not prove that the two values coincide, nor does it resolve the Open Question $11$ for the broader class $\mathcal{F}_p$.
\end{rem}  
\end{enumerate}

\bigskip

\noindent \textbf{Acknowledgements:} I would like to thank Carlo Mantegazza for the discussions we’ve had about this project and for his constant support over the years.


\end{document}